\documentclass[12pt]{article}

\usepackage{graphicx,psfrag}
\usepackage{psfrag,eepic,epsfig}
\usepackage{sectsty}
\usepackage{setspace}
\usepackage{color}

\usepackage{setspace}
\usepackage{amsmath, epsfig, cite, ifpdf}
\usepackage{amssymb}
\usepackage{amsfonts}
\usepackage{latexsym}
\usepackage{graphicx,psfrag,eepic,rawfonts}
\usepackage{enumerate}
\usepackage{indentfirst}
\input{epsf}

\newtheorem{thm}{Theorem}[section]

\newtheorem{cor}[thm]{Corollary}

\newtheorem{lem}[thm]{Lemma}
\newtheorem{conj}[thm]{Conjecture}
\newtheorem{problem}[thm]{Problem}
\newtheorem{exam}[thm]{Example}

\newcommand{\qed}{{\hfill\rule{4pt}{7pt}}}
\def\pf{\noindent {\it Proof.} }

\numberwithin{equation}{section}

\makeatletter \@addtoreset{equation}{section} \makeatother

\title {\bf More on the skew-spectra of bipartite graphs and Cartesian products of
graphs\footnote{Supported by NSFC and the ``973" program. }}

\author{
{\small Xiaolin Chen, Xueliang Li, Huishu Lian}\\
{\small Center for Combinatorics and LPMC-TJKLC}\\
{\small Nankai University, Tianjin 300071, P.R. China}\\
{\small E-mail: chxlnk@163.com; lxl@nankai.edu.cn; lhs6803@126.com}
   }
\date{}
\begin{document}

\maketitle

\begin{abstract}
Given a graph $G$, let $G^\sigma$ be an oriented graph of $G$ with
the orientation $\sigma$ and skew-adjacency matrix $S(G^\sigma)$.
Then the spectrum of $S(G^\sigma)$ is called the skew-spectrum of
$G^\sigma$, denoted by $Sp_S(G^\sigma)$. It is known that a graph
$G$ is bipartite if and only if there is an orientation $\sigma$ of
$G$ such that $Sp_S(G^\sigma)=iSp(G)$. In [D. Cui, Y. Hou, On the
skew spectra of Cartesian products of graphs, Electron. J. Combin.
20(2013), \#P19], Cui and Hou conjectured that such orientation of a
bipartite graph is unique under switching-equivalence. In this
paper, we prove that the conjecture is true. Moreover, we give an
orientation of the Cartesian product of a bipartite graph and a
graph, and then determine the skew-spectrum of the resulting
oriented product graph, which generalizes Cui and Hou's result, and
can be used to construct more oriented graphs with maximum skew energy. \\


\end{abstract}

\section{Introduction}

Let $G$ be a simple undirected graph on order $n$ with vertex set
$V(G)$ and edge set $E(G)$. Suppose $V(G)=\{v_1,v_2,\ldots,v_n\}$.
Then the {\it adjacency matrix} of $G$ is the $n\times n$ symmetric
matrix $A(G)=[a_{ij}]$, where $a_{ij}=1$ if the vertices $v_i$ and
$v_j$ are adjacent, and $a_{ij}=0$ otherwise. The {\it spectrum} of
$G$, denoted by $Sp(G)$, is defined as the spectrum of $A(G)$.

Let $G^\sigma$ be an oriented graph of $G$ with the orientation
$\sigma$, which assigns to each edge of $G$ a direction so that the
induced graph $G^\sigma$ becomes an oriented graph or a directed
graph. Then $G$ is called the underlying graph of $G^\sigma$. The
{\it skew-adjacency matrix} of $G^\sigma$ is the $n\times n$ matrix
$S(G^\sigma)=[s_{ij}]$, where $s_{ij}=1$ and $s_{ji}=-1$ if $\langle
v_i,v_j\rangle$ is an arc of $G^\sigma$, otherwise
$s_{ij}=s_{ji}=0$. Obviously, $S(G^\sigma)$ is a skew-symmetric
matrix, and thus all the eigenvalues are purely imaginary numbers or
0, which form the spectrum of $S(G^\sigma)$ and are said to be the
{\it skew-spectrum} of $G^\sigma$. The eigenvalues of $S(G^\sigma)$
are called the {\it skew eigenvalues} of $G^\sigma$.

The {\it energy} $\mathcal{E}(G)$ of an undirected graph is defined
as the sum of the absolute values of all the eigenvalues of $G$,
which was introduced by Gutman in \cite{G}. We refer the survey
\cite{GLZ} and the book \cite{LSG} to the reader for details. The
{\it skew energy} of an oriented graph, as one of various
generalizations of the graph energy, was first proposed by Adiga et
al. \cite{ABC}. It is defined as the sum of the absolute values of
all the eigenvalues of $S(G^\sigma)$, denoted by
$\mathcal{E}_S(G^\sigma)$. Most of the results on the skew energy
are collected in our recent survey \cite{LL}.

In \cite{ABC}, Adiga et al. obtained some properties of the skew energy
and proposed some open problems, such as the following two problems.

\begin{problem}\label{prob1}
Find new families of oriented graphs $G^\sigma$ with $\mathcal{E}_S(G^\sigma)=\mathcal{E}(G)$.
\end{problem}
\begin{problem}\label{prob2}
Which $k$-regular graphs on $n$ vertices have orientations $G^\sigma$ with
$\mathcal{E}_S(G^\sigma)= n\sqrt{k}$, or equivalently, $S(G^\sigma)^{T}S(G^\sigma)=k I_n$ ?
\end{problem}

For Problem \ref{prob1}, it is clear that if an oriented graph
$G^\sigma$ satisfies $Sp_S(G^\sigma)=iSp(G)$ then
$\mathcal{E}_S(G^\sigma)=\mathcal{E}(G)$. In \cite{HL,Shader}, they
proved that $Sp_S(G^\sigma)=iSp(G)$ for any orientation $\sigma$ of
$G$ if and only if $G$ is a tree, and $Sp_S(G^\sigma)=iSp(G)$ for
some orientation $\sigma$ of $G$ if and only if $G$ is a bipartite
graph. They also pointed out that the {\it elementary orientation}
of a bipartite graph $G=G(X,Y)$, which assigns each edge the
direction from $X$ to $Y$,  is such an orientation that
$Sp_S(G^\sigma)=iSp(G)$. Are there any other orientations ?

Let $W$ be a vertex subset of an oriented graph $G^\sigma$ and
$\overline{W}=V(G^\sigma)\backslash W$. Another oriented graph
$G^{\sigma'}$ of $G$, obtained from $G^\sigma$ by reversing the
orientations of all arcs between $W$ and $\overline{W}$, is said to
be obtained from $G^\sigma$ by {\it switching with respect to $W$}.
Two oriented graphs $G^\sigma$ and $G^{\sigma'}$ are said to be
switching-equivalent if $G^{\sigma'}$ can be obtained from
$G^\sigma$ by a sequence of switchings. Note that if $G^\sigma$ and
$G^{\sigma'}$ are {\it switching-equivalent}, then
$\mathcal{E}_S(G^\sigma)=\mathcal{E}_S(G^{\sigma'})$; see
\cite{HSZ}. Moreover, in a recent paper \cite{CH}, Cui and Hou posed
the following conjecture:
\begin{conj}\label{conj}
Let $G=G(X,Y)$ be a bipartite graph and $\sigma$ be an orientation of $G$.
Then $Sp_S(G^\sigma)=iSp(G)$ if and only if $\sigma$ is switching-equivalent
to the elementary orientation of $G$.
\end{conj}

In Section \ref{bipartite}, we prove that the Conjecture \ref{conj}
is true. We also obtain some other results of the skew-spectrum of
an oriented bipartite graph.

For Problem \ref{prob2}, it is known that
$\mathcal{E}_S(G^\sigma)\leq n\sqrt{\Delta}$ and equality holds if
and only if $S(G^\sigma)^TS(G^\sigma)=\Delta I_n$, which implies
that $G$ is regular. Some families of oriented graphs with maximum
skew energy were characterized in \cite{ABC,GX,Tian}. Moreover, in
\cite{CH} Cui and Hou constructed a new family of oriented graphs
with maximum skew energy by considering the skew-spectrum of
Cartesian product $P_2\Box G$. In Section \ref{Cartesian}, we extend
their result of the product graph $P_2\Box G$ to that of the product
graph $H\Box G$, where $H$ is a bipartite graph. Using it, we obtain
a larger new family of oriented graphs with maximum skew energy.

\section{Oriented bipartite graphs with $Sp_S(G^\sigma)=iSp(G)$ }\label{bipartite}

In this section, we consider bipartite graphs and their orientations.
We prove that Conjecture \ref{conj} is true and also obtain the characterizations
of orientations of a bipartite graph $G$ with $Sp_S(G^\sigma)=iSp(G)$.

First, let's recall some definitions. Let $G^\sigma$ be an oriented
graph and $C_{2\ell}$ be an undirected even cycle of $G$. Then
$C_{2\ell}$ is said to be {\it evenly oriented} relative to
$G^\sigma$ if it has an even number of edges oriented in clockwise
direction (and now it also has an even number of edges oriented in
anticlockwise direction, since $C_{2\ell}$ is an even cycle);
otherwise $C_{2\ell}$ is called {\it oddly oriented}. The even cycle
$C_{2\ell}$ is said to be  {\it oriented uniformly} if $C_{2\ell}$
is oddly (resp., evenly) oriented relative to $G^\sigma$ when $\ell$
is odd (resp., even). It should be noted that if an even cycle
$C_{2\ell}$ is oriented uniformly in $G^\sigma$, then after
switching operations on $G^\sigma$, $C_{2\ell}$ is also oriented
uniformly. The following lemma was obtained in \cite{CH}.
\begin{lem}\cite{CH} \label{even}
Let $G$ be a bipartite graph and $\sigma$ be an orientation of $G$.
Then $Sp_S(G^\sigma)=iSp(G)$ if and only if every even cycle is
oriented uniformly in $G^\sigma$.
\end{lem}

Now we prove the following result which implies that Conjecture
\ref{conj} is true.
\begin{thm}
Let $G=G(X,Y)$ be a bipartite graph and $\sigma$ be an orientation
of $G$. Then $Sp_S(G^\sigma)=iSp(G)$ if and only if $\sigma$ is
switching-equivalent to the elementary orientation of $G$.
\end{thm}
\pf By Lemma \ref{even}, we can easily get the sufficiency of the
theorem. We prove the necessity by induction on the number $m$ of
edges in $G$.

When $m=1$, the case is trivial. Assume that the result holds for
$m-1$. Let $G=G(X,Y)$ be a bipartite graph with $m$ edges and
$\sigma$ be an orientation of $G$ such that $Sp_S(G^\sigma)=iSp(G)$.
Then by Lemma \ref{even}, every even cycle of $G^\sigma$ is oriented
uniformly. Suppose that $e=xy$ is an edge of $G$ with $x\in X$,
$y\in Y$ and $\hat{e}$ is the corresponding arc of $G^\sigma$.
Consider the oriented bipartite graph $G^\sigma-\hat{e}$. Note that
every even cycle of $G^\sigma-\hat{e}$ is also oriented uniformly,
and thus $Sp_S(G^\sigma-\hat{e})=iSp(G-e)$. By induction, we obtain
that the orientation of $G^\sigma-\hat{e}$ is switching-equivalent
to the elementary orientation of $G-e$. That is, there exists a
sequence of switchings which transform the orientation of
$G^\sigma-\hat{e}$ to the elementary orientation of $G-e$. Applying
the same switching operations on $G^\sigma$, we obtain an oriented
graph $G^{\sigma'}$, whose arcs except the arc corresponding to $e$
have directions from $X$ to $Y$.

If the direction of the arc corresponding to $e$ in $G^{\sigma'}$ is from $x$ to $y$,
it follows that $\sigma$ is switching-equivalent to the elementary orientation of $G$.

If the direction of the arc corresponding to $e$ in $G^{\sigma'}$ is
from $y$ to $x$. We claim that $e$ is a cut edge of $G$. Otherwise,
there exists a cycle containing $e$, say $C_{2k}=x_0y_0x_1y_1\cdots
x_{k-1}y_{k-1}x_0$ in clockwise direction. We find that $C_{2k}$
contains precisely $k-1$ arcs in clockwise direction, which
contradicts that $C_{2k}$ is oriented uniformly. Hence $e$ is a cut
edge of $G$. Then $V(G)$ can be partitioned into two parts $W_1$ and
$W_2$ such that $e$ is the only edge between them. By switching with
respect to $W_1$ in $G^{\sigma'}$, the direction of the arc
corresponding to $e$ is reversed and the directions of other arcs
keep unchanged. Then we conclude that $\sigma$ is
switching-equivalent to the elementary orientation of $G$. The proof
is thus complete.\qed

Lemma \ref{even} provides a good characterization of the oriented
bipartite graphs $G^\sigma$ with $Sp_S(G^\sigma)=iSp(G)$. But it
requires one to check that every even cycle is oriented uniformly. A
natural question is how to simplify the checking task, such as only
to check all chordless cycles. Based on this, we get the following
result.
\begin{thm}\label{chordless}
Let $G$ be a bipartite graph and $\sigma$ be an orientation of $G$. If every chordless
cycle of $G^\sigma$ is oriented uniformly, then $Sp_S(G^\sigma)=iSp(G)$.
\end{thm}
\pf We prove the theorem by contradiction. Let $G$ be a bipartite
graph and $\sigma$ be an orientation of $G$ such that every
chordless cycle is oriented uniformly. But $Sp_S(G^\sigma)\neq
iSp(G)$. By Lemma \ref{even}, there exists an even cycle that is not
oriented uniformly. Choose a shortest cycle $C_{2\ell}$ that is not
oriented uniformly, that is, $C_{2\ell}$ is oddly oriented in
$G^\sigma$ if $\ell$ is even and evenly oriented if $\ell$ is odd.
It follows that $C_{2\ell}$ contains a chord $x_1y_1$. Suppose that
$C_{2\ell}=x_1x_2\cdots x_{\ell_1}y_1y_2\cdots y_{2\ell-\ell_1}x_1$
in clockwise direction. Consider the two cycles $C_1=x_1x_2\cdots
x_{\ell_1}y_1x_1$ with length $\ell_1+1$ in clockwise direction and
$C_2=x_1y_1y_2\cdots y_{2\ell-\ell_1}x_1$ with length
$2\ell-\ell_1+1$ in clockwise direction. Note that $C_1$ and $C_2$
must be even. Suppose that $C_1$ and $C_2$ contain $r_1$ and $r_2$
arcs in clockwise direction, respectively. Since $C_1$ and $C_2$ are
oriented uniformly, we get
\begin{equation*}
\frac{\ell_1+1}{2}\equiv r_1 \mod 2,
\hspace{10pt}\text{and}\hspace{10pt} \frac{2\ell-\ell_1+1}{2}\equiv
r_2 \mod 2.
\end{equation*}
It follows that $\ell+1\equiv r_1+r_2\mod 2$. Observe that if the
arc corresponding to $x_1y_1$ is clockwise in $C_1$, then it must be
anticlockwise in $C_2$, and vice versa. Thus $C_{2\ell}$ contains
$r_1+r_2-1$ arcs in clockwise direction. Therefore, $C_{2\ell}$ is
oriented uniformly, a contradiction. The proof is now complete.\qed

Combining with Lemma \ref{even}, we immediately obtain the following corollary.
\begin{cor}
Let $G$ be a bipartite graph and $\sigma$ be an orientation of $G$. Then $Sp_S(G^\sigma)=iSp(G)$
if and only if all chordless cycles are oriented uniformly in $G^\sigma$.
\end{cor}

\noindent {\bf Remark 2.5} Let $\mathcal{C}$ denote the set of all
cycles of a bipartite graph $G$. A subset $\mathcal{S}$ of
$\mathcal{C}$ is called a {\it generating set} of $\mathcal{C}$ if
for any cycle $C$ of $\mathcal{C}$, $C\in \mathcal{S}$ or there is a
sequence of cycles $C_1, C_2, \ldots, C_k$ in $\mathcal{S}$ such
that $C=((C_1\vartriangle C_2)\vartriangle C_3)\cdots \vartriangle
C_k$ and $C_1\vartriangle C_2$, $(C_1\vartriangle C_2)\vartriangle
C_3$, $\ldots$, $((C_1\vartriangle C_2)\vartriangle C_3)\cdots
\vartriangle C_{k-1}$ all are cycles. With this notation, one can
prove that for an oriented bipartite graph $G^\sigma$,
$Sp_S(G^\sigma)=iSp(G)$ if and only if every cycle in a generating
set $\mathcal{S}$ of $\mathcal{C}$ is oriented uniformly in
$G^\sigma$. Actually, the set of chordless cycles of a graph is a
generating set of the set of all cycles of $G$.

\section{The skew-spectrum of $H\Box G$ with $H$ bipartite}\label{Cartesian}

In this section, we give an orientation of the Cartesian product
$H\Box G$, where $H$ is bipartite, by extending the orientation of
$P_m\Box G$ in \cite{CH}, and we calculate its skew-spectrum. As an
application of this orientation, we construct a larger new family of
oriented graphs with maximum skew energy, which generalizes the
construction in \cite{CH}.

Let $H$ and $G$ be graphs with $m$ and $n$ vertices, respectively.
The Cartesian product $H\square G$ of $H$ and $G$ is a graph with
vertex set $V(H)\times V(G)$ and there exists an edge between $(u_1,
v_1)$ and $(u_2, v_2)$ if and only if $u_1=u_2$ and $v_1 v_2$ is an
edge of $G$, or $v_1=v_2$ and $u_1 u_2$ is an edge of $H$. Assume
that $H^\tau$ is any orientation of $H$ and $G^\sigma$ is any
orientation of $G$.  There is a natural way to give an orientation
$H^\tau \square G^\sigma$ of $H^\tau$ and $G^\sigma$. There is an
arc from $(u_1, v_1)$ to $(u_2, v_2)$ if and only if $u_1=u_2$ and
$(v_1, v_2)$ is an arc of $G^\sigma$, or $v_1=v_2$ and $(u_1, u_2)$
is an arc of $H^\tau$.

When $H$ is a bipartite graph with bipartition $X$ and $Y$, we
modify the above orientation of $H^\tau\square G^\sigma$ with the
following method. If there is an arc from $(u, v_1)$ to $(u, v_2)$
in $H^\tau\square G^\sigma$ and $u\in Y$, then we reverse the
direction of the arc. The other arcs keep unchanged. This new
orientation of $H\square G$ is denoted by $(H^\tau\square
G^\sigma)^o$.

\begin{thm}\label{spectrum}
Let $H^\tau$ be an oriented bipartite graph of order $m$ and let the
skew eigenvalues of $H^\tau$ be the non-zero values $\pm \mu_1i,\pm
\mu_2i,\dots,\pm \mu_ti$ and $m-2t$ $0$'s. Let $G^\sigma$ be an
oriented graph of order $n$ and let the skew eigenvalues of
$G^\sigma$ be the non-zero values $\pm \lambda_1i,\pm
\lambda_2i,\dots,\pm \lambda_ri$ and $n-2r$ $0$'s. Then the skew
eigenvalues of the oriented graph $(H^\tau\Box G^\sigma)^o$ are $\pm
i\sqrt{\mu_j^2+\lambda_k^2}$ with multiplicities $2$, $j=1,\dots,t$,
$k=1,\dots,r$, $\pm\mu_ji$ with multiplicities $n-2r$,
$j=1,\dots,t$, $\pm\lambda_ki$ with multiplicities $m-2t$,
$k=1,\dots,r$, and $0$ with multiplicities $(m-2t)(n-2r)$.
\end{thm}
\pf Let $H=H(X,Y)$ be a bipartite graph with $|X|=m_1$ and
$|Y|=m_2$. With suitable labeling of the vertices of $H\Box G$, the
skew-adjacency matrix $S=S((H^\tau\Box G^\sigma)^o)$ can be
formulated as follows:
$$S=I'_{m_1+m_2}\otimes
S(G^\sigma)+S(H^\tau)\otimes I_n,$$ where $I'_{m_1+m_2}=(a_{ij})$,
$a_{ii}=1$ if $1\leq i\leq m_1$, $a_{ii}=-1$ if $m_1+1\leq i\leq m$,
and $a_{ij}=0$ for $i\neq j$; $S(H^\tau)$ is the partition matrix
$\left(\begin{array}{cc} 0& B\\ -B^T& 0\end{array}\right)$ and $B$
is an $m_1\times m_2$ matrix.

We first determine the singular values of $S$. Note that $S$,
$S(H^\tau)$ and $S(G^\sigma)$ are all skew symmetric. By
calculation, we have
\begin{equation*}
\begin{split}
SS^T=&\left(I'_{m_1+m_2}\otimes S(G^\sigma)+S(H^\tau)\otimes I_n\right)(I'_{m_1+m_2}\otimes (-S(G^\sigma))+(-S(H^\tau))\otimes I_n)\\
=&-\left[\left(I_{m_1+m_2}\otimes S^2(G^\sigma)+S^2(H^\tau)\otimes I_n\right)+(I'_{m_1+m_2}\otimes S(G^\sigma))(S(H^\tau)\otimes I_n)\right.\\
&\left.+(S(H^\tau)\otimes I_n)(I'_{m_1+m_2}\otimes
S(G^\sigma))\right].
\end{split}
\end{equation*}
Define $\sigma_i=1$ for $i=1,2,\dots,m_1$ and $\sigma_i=-1$ for
$i=m_1+1,m_1+2,\dots,m$. Denote $M^1=[M_{ij}^1]=(I'_{m_1+m_2}\otimes
S(G^\sigma))(S(H^\tau)\otimes I_n)$ and
$M^2=[M_{ij}^2]=(S(H^\tau)\otimes I_n)(I'_{m_1+m_2}\otimes
S(G^\sigma))$. Note that $M^1$ and $M^2$ are both $m\times m$
partition matrix in which every entry is an $n\times n$ submatrix.
Direct computing gives
$$M_{ij}^1+M_{ij}^2=S(H^\tau)_{ij}S(G^\sigma)((-1)^{\sigma_i}+(-1)^{\sigma_j}).$$

For any $1\leq i,j\leq m$, if $S(H^\tau)_{ij}=0$, then
$M_{ij}^1+M_{ij}^2=0$. Otherwise the vertices corresponding to $i$
and $j$ in $H^\tau$ are in different parts of the bipartition. That
is, $1\leq i\leq m_1, m_1+1\leq j\leq m$ or $1\leq j\leq m_1,
m_1+1\leq i\leq m$. Then $(-1)^{\sigma_i}+(-1)^{\sigma_j}=0$, an
thus $M_{ij}^1+M_{ij}^2=0$. It follows that $M^1+M^2=0$. Hence,
$$SS^T=-\left(I_{m_1+m_2}\otimes S^2(G^\sigma)+S^2(H^\tau)\otimes I_n\right).$$
Therefore, the eigenvalues of $SS^T$ are
$\mu(H^\tau)^2+\lambda(G^\sigma)^2$, where $\mu(H^\tau)i\in
Sp_S(H^\tau)$ and $\lambda(G^\sigma)i\in Sp_S(G^\sigma)$. Then the
skew-spectrum of $(H^\tau\Box G^\sigma)^o$ follows. The proof is
thus complete. \qed

As an application of Theorem \ref{spectrum}, we can now construct a
new family of oriented graphs with maximum skew energy.

\begin{thm} \label{spectrum1}
Let $H^\tau$ be an oriented $\ell$-regular bipartite graph on $m$
vertices with maximum skew energy $\mathcal{E_S}(H^\tau)=m
\sqrt{\ell}$ and $G^\sigma$ be an oriented $k$-regular graph on $n$
vertices with maximum skew energy $\mathcal{E_S}(G^\sigma)=n
\sqrt{k}$. Then the oriented graph $(H^\tau\Box G^\sigma)^o$ of
$H\Box G$ has the maximum skew energy $\mathcal{E_S}((H^\tau\Box
G^\sigma)^o)=mn \sqrt{\ell+k}$.
\end{thm}
\pf Since $H^\tau$ and $G^\sigma$ have maximum skew energy,
$S(H^\tau)S(H^\tau)^T=\ell I_m$ and
$S(G^\sigma)S(G^\sigma)^T$\\$=kI_n$. Then the skew eigenvalues of
$H^\tau$ are all $\pm i\sqrt{\ell}$ and the skew eigenvalues of
$G^\sigma$ are all $\pm i\sqrt{k}$. By Theorem \ref{spectrum},
$(H^\tau\Box G^\sigma)^o$ have all skew eigenvalues $\pm
imn\sqrt{\ell+k}$. \qed

The following result was obtained in \cite{CH}, which can be viewed
as an immediate corollary of Theorem \ref{spectrum1}.
\begin{cor} \cite{CH}
Let $G^\sigma$ be an oriented $k$-regular graph on $n$ vertices with
maximum skew energy $\mathcal{E_S}(G^\sigma)=n \sqrt{k}$. Then the
oriented graph $(P_2 \Box G^\sigma)^o$ of $P_2\Box G$ has maximum
skew energy $\mathcal{E_S}((P_2\Box G^\sigma)^o)=2n \sqrt{k+1}$.
\end{cor}

Adiga et al. \cite{ABC} showed that a $1$-regular connected graph
that has an orientation with maximum skew energy is $K_2$; while a
$2$-regular connected graph has an orientation with maximum skew
energy if and only if it is $C_4$ with oddly orientation. Tian
\cite{Tian} proved that there exists a $k$-regular graph with
$n=2^k$ vertices having an orientation $\sigma$ with maximum skew
energy. Cui and Hou \cite{CH} constructed a $k$-regular graph of
order $n=2^{k-1}$ having an orientation $\sigma$ with maximum skew
energy. The following examples provide new families of oriented
graphs with maximum skew energy that have much less vertices.
\begin{exam}
Let $G_1=K_{4,4}$, $G_2=K_{4,4}\Box G_1, \dots, G_r=K_{4,4}\Box
G_{r-1}$. Because there is an orientation of $K_{4,4}$ with maximum
skew energy 16; see \cite{CLL}. Thus, we can get an orientation of
$G_r$ with maximum energy $2^{3r}\sqrt{4r}$. This provides a family
of $4r$-regular graphs of order $n=2^{3r}$ having an orientation
with skew energy $2^{3r}\sqrt{4r}$ for $r\geq 1$.
\end{exam}
\begin{exam}
Let $G_1=K_4$, $G_2=K_{4,4}\Box G_1,\dots,G_r=K_{4,4}\Box G_{r-1}$.
It is known that $K_4$ has an orientation with maximum skew energy;
see \cite{GX}. Thus we can get an orientation of $G_r$ with maximum
energy $2^{3r-1}\sqrt{4r-1}$. This provides a family of
$4r-1$-regular graphs of order $n=2^{3r-1}$ having an orientation
with skew energy $2^{3r-1}\sqrt{4r-1}$ for $r\geq 1$.
\end{exam}
\begin{exam}
Let $G_1=C_4$, $G_2=K_{4,4}\Box G_1, \dots, G_r=K_{4,4}\Box
G_{r-1}$. Thus, we can get an orientation of $G_r$ with maximum
energy $2^{3r-1}\sqrt{4r-2}$. This provides a family of
$4r-2$-regular graphs of order $n=2^{3r-1}$ having an orientation
with skew energy $2^{3r-1}\sqrt{4r-2}$ for $r\geq 1$.
\end{exam}
\begin{exam}
Let $G_1=P_2$, $G_2=K_{4,4}\Box G_1, \dots, G_r=K_{4,4}\Box
G_{r-1}$. Thus, we can get an orientation of $G_r$ with maximum
energy $2^{3r-2}\sqrt{4r-3}$. This provides a family of
$4r-3$-regular graphs of order $n=2^{3r-2}$ having an orientation
with skew energy $2^{3r-2}\sqrt{4r-3}$ for $r\geq 1$.
\end{exam}

\end{document}